# NUMERICAL APPROXIMATIONS OF THE CAHN-HILLIARD AND ALLEN-CAHN EQUATIONS WITH GENERAL NONLINEAR POTENTIAL USING THE INVARIANT ENERGY QUADRATIZATION APPROACH

XIAOFENG YANG* AND GUO-DONG ZHANG†

**Abstract.** In this paper, we carry out stability and error analyses for two first-order, semi-discrete time stepping schemes, which are based on the newly developed Invariant Energy Quadratization approach, for solving the well-known Cahn-Hilliard and Allen-Cahn equations with general nonlinear bulk potentials. Some reasonable sufficient conditions about boundedness and continuity of the nonlinear functional are given in order to obtain optimal error estimates. These conditions are naturally satisfied by two commonly used nonlinear potentials including the double-well potential and regularized logarithmic Flory-Huggins potential. The well-posedness, unconditional energy stabilities and optimal error estimates of the numerical schemes are proved rigorously.

**1. Introduction.** In this paper, we carry out stability analyses and error estimates for two first-order, semi-discrete time-stepping numerical schemes, that are based on the newly developed Invariant Energy Quadratization approach, for solving the well-known Cahn-Hilliard and Allen-Cahn equations with general nonlinear bulk potentials. These two equations, produced by the gradient flow (a diffusive process) in $L^2$ and $H^{-1}$ spaces respectively, are typical representatives of the phase field (diffusive interface) method, which is a well-known efficient modeling and a numerical tool to resolve the motion of free interfaces between multiple material components. A main advantage of this approach is that it permits to solve the free interface problem by integrating a set of partial differential equations, thus avoiding the explicit treatment of the boundary conditions at the interface, see [3, 4, 10, 19, 20, 23, 25, 27, 32] and the references therein about the theoretical/numerical work, as well as their wide applications in various science and engineering fields.

For algorithms design of any phase field models related to the Allen-Cahn and Cahn-Hilliard equations, a significant goal is to verify the energy stable property at the discrete level irrespectively of the coarseness of the discretization in time and space. In what follows, those algorithms will be called unconditionally energy stable. Schemes with this property are especially preferred for solving phase field models since it is not only critical for the numerical scheme to capture the correct long-time dynamics of the system, but also provides sufficient flexibility for dealing with the stiffness issue, induced by the thin interface. Meanwhile, since the dynamics of coarse-graining (macroscopic) process may undergo rapid changes near the interface, the noncompliance of energy dissipation laws may lead to spurious numerical solutions if the mesh and time step size are not carefully controlled.

However, it is a challenging task to develop unconditionally energy stable schemes to resolve the stiffness issue from the thin interface. It is shown that the traditional fully-implicit or explicit discretization for the nonlinear term can cause very severe time step constraint (called conditionally energy stable) on the interfacial width [1, 13, 31]. Many efforts had been done (cf. [5, 12–14, 18, 24, 26, 28, 29, 31, 33, 34, 34, 38, 45, 47, 50] and the references therein) in order to remove this constraint. Among those, the *convex splitting approach* [12, 29, 33, 35, 46] and *linear stabilization approach* [5, 24, 31, 34, 38, 45, 50] are the two most commonly used numerical techniques. In the convex splitting approach, the convex part of the potential is treated implicitly and the concave part is treated explicitly. The scheme is unconditionally energy stable, however, it usually produces nonlinear type schemes, thus the implementation is complicated and the computational cost might be high. The linear stabilization approach treats the nonlinear term explicitly thus it is efficient and very easy to implement. In order to remove the time step constraint dependence on the interfacial width, a linear stabilizing term is added and its magnitude usually depends on the upper bound of the second order derivative of the nonlinear potential. Unfortunately, the upper bound is usually infinite, therefore the remedy is to modify the nonlinear potential, that in turn results in additional accuracy issues. Moreover, it appears difficult to design second-order unconditionally energy stable schemes with a stabilized approach, although some progress has been made recently in [21, 22].

---

*Corresponding author, Department of Mathematics, University of South Carolina, Columbia, SC, 29208, USA. Email: xfyang@math.sc.edu. This author's research is partially supported by the U.S. National Science Foundation under grant numbers DMS-1418898 and DMS-1720212.

†Department of Mathematics, Yantai University, Yantai, 264005, Shandong, P. R. China. Email: gdzhang@ytu.edu.cn. This author's research is partially supported by National Science Foundation of China under grant numbers 11601468 and 11771375.

Recently, a novel numerical method, called *Invariant Energy Quadratization* (IEQ) approach, has been developed and successfully applied to solve a variety of gradient flow models (cf. [6, 8, 15, 17, 36, 37, 39–44, 44, 48, 49]). Its essential idea is to transform the bulk potential into a quadratic form (since the nonlinear potential is usually bounded from below) using a set of new variables via the change of variables. Then, for the reformulated model in terms of the new variables that still retains the identical energy dissipation law, all nonlinear terms can be treated semi-explicitly, which then yields a linear system. This method bypasses those typical challenges such as the justification/adjustment of convexity or implicit/explicit terms, to arrive at first, second order or even higher order in time unconditionally energy stable schemes readily.

Although one might think that it would be natural to derive the corresponding error analysis for the IEQ type schemes by analogy with the proof of stability, the reality is quite the opposite. An exceptional case is for the double well potential where many terms can be simplified (cf. Remark 3.2). For general nonlinear potentials, an essential difficulty comes from the way of quadratization to introduce the new variable, that actually leads the new variable to act as a *encapsulation*, making it difficult to estimate the quantitative relation between the new and original variables. To the best of the authors' knowledge, we are not aware of any results about the error analysis of IEQ type schemes with general nonlinear potentials and almost all works had been focused on its remarkable unconditional energy stability. In view of the scarce of error analysis, the main objective of this paper is to derive optimal error estimates for the IEQ schemes for solving the Cahn-Hilliard and Allen-Cahn equations. We give some reasonable sufficient conditions about boundedness and continuity for the nonlinear functionals in order to obtain optimal error estimates. These conditions are naturally satisfied by the commonly used double well potential and regularized logarithmic Flory-Huggins potential. Moreover, the analytical approach developed in this paper is general enough and thus it can work as a standard framework to derive error estimates of IEQ type schemes for various gradient flow models with diverse nonlinear potentials.

The rest of paper is organized as follows. In Section 2, we give a brief introduction to the Cahn-Hilliard and Allen-Cahn systems. In Section 3, for the IEQ scheme of solving the fourth order Cahn-Hilliard equation, we study the unconditional energy stability, practical implementation, wellposedness, and derive the optimal error estimates. In Section 4, similar analytical work is performed for the second order Allen-Cahn equation. In Section 5, some concluding remarks are given.

**2. PDE Models and their energy laws.** We consider the following Lyapunov energy functional,

$$E(\phi) = \int_\Omega \left(\frac{\epsilon^2}{2}|\nabla\phi|^2 + F(\phi)\right) d\boldsymbol{x}, \tag{2.1}$$

where $\phi(\boldsymbol{x}, t)$ is the unknown scalar function, $\boldsymbol{x} \in \Omega \subseteq \mathbb{R}^d$ ($d = 2, 3$), $F(\phi)$ is the nonlinear bulk potential, $\epsilon$ is an interface/penalty parameter causing *stiffness* issue into the PDE system when $\epsilon \ll 1$. There are two commonly used nonlinear bulk potentials for $F(\phi)$:

(i) Ginzburg-Landau double-well type potential, cf. [3, 23, 28]:

$$F_{db}(x) = \frac{1}{4}(x^2 - 1)^2, x \in (-\infty, \infty), \tag{2.2}$$

(ii) Logarithmic Flory-Huggins potential, cf. [2, 3, 9, 11]:

$$F_{fh}(x) = x\ln x + (1-x)\ln(1-x) + \theta(x - x^2), \ \theta > 0, x \in (0, 1). \tag{2.3}$$

For either of these two nonlinear potentials, we note there always exists a positive constant $A$ such

that

$$\begin{cases} F_{db}(x) > -A, & \forall x \in (-\infty, \infty); \\ F_{fh}(x) > -A, & \forall x \in (0,1), \end{cases} \tag{2.4}$$

where we can simply choose $A = 1$ for both cases.

By applying the varational approach for the total free energy (2.1) in $H^{-1}(\Omega)$, we obtain the Cahn-Hilliard type system that reads as:

$$\phi_t - M\Delta w = 0, \tag{2.5}$$
$$w = -\epsilon^2 \Delta \phi + f(\phi), \quad (\boldsymbol{x}, t) \in \Omega \times (0, T], \tag{2.6}$$

where $M$ is the mobility constant, $w$ is the chemical potential, and $f(\phi) = F'(\phi)$. The initial condition is $\phi|_{(t=0)} = \phi_0$. For simplicity, we choose the suitable boundary conditions so that all complexities from the boundary integrals can be removed, i.e.,

$$(i) \ \phi, w \ are \ periodic; \ or \ (ii) \ \partial_{\boldsymbol{n}}\phi|_{\partial\Omega} = \partial_{\boldsymbol{n}}w|_{\partial\Omega} = 0. \tag{2.7}$$

By applying the varational approach for the total free energy (2.1) in $L^2(\Omega)$, we obtain the Allen-Cahn type system that reads as:

$$\phi_t + M(-\epsilon^2 \Delta \phi + f(\phi)) = 0, \quad (\boldsymbol{x}, t) \in \Omega \times (0, T]. \tag{2.8}$$

Its boundary conditions are

$$(i) \ \phi \ is \ periodic; \ or \ (ii) \ \partial_{\boldsymbol{n}}\phi|_{\partial\Omega} = 0. \tag{2.9}$$

An important feature of the Cahn-Hilliard and Allen-Cahn equations is that they both satisfy energy dissipation law. For the Cahn-Hilliard system (2.5)-(2.6), by taking the $L^2$ inner product of (2.5) with $-w$, of (2.6) with $\phi_t$, performing integration by parts and summing up two equalities, we obtain

$$\frac{d}{dt}E(\phi) = -M\|\nabla w\|^2 \leq 0. \tag{2.10}$$

For the Allen-Cahn system (2.8), by taking the $L^2$ inner product of with $\phi_t$, performing integration by parts, we obtain

$$\frac{d}{dt}E(\phi) = -\frac{1}{M}\|\phi_t\|^2 \leq 0. \tag{2.11}$$

**3. Cahn-Hilliard equation.** We first introduce some notations that will be used throughout the paper. We let $L^p(\Omega)$ denote the usual Lebesgue space on $\Omega$ with the norm $\|\cdot\|_{L^p}$. The inner product and norm in $L^2(\Omega)$ are denoted by $(\cdot, \cdot)$ and $\|\cdot\|$, respectively. $W^{k,p}(\Omega)$ stands for the standard Sobolev spaces equipped with the standard Sobolev norms $\|\phi\|_{k,p}$. For $p = 2$, we write $H^k(\Omega)$ for $W^{k,2}(\Omega)$, and the corresponding norm is $\|\phi\|_k$. We define two Sobolev spaces:

$$H_{per}(\Omega) = \{\phi \ is \ periodic, \phi \in H^1(\Omega) \ and \ \int_\Omega \phi d\boldsymbol{x} = 0\},$$

$$H(\Omega) = \{\phi \in H^1(\Omega) \ and \ \int_\Omega \phi d\boldsymbol{x} = 0\}.$$

**3.1. Unconditional energy stable linear scheme using the IEQ approach.** We recall that the main challenge to develop efficient, unconditionally energy stable schemes for solving the system (2.5)-(2.6) lies in how to discretize the nonlinear term $f(\phi)$. We now consider the newly developed IEQ approach, see [6, 8, 15, 17, 36, 37, 39–44, 44, 48, 49]. Note $F(\phi)$ is bounded from below as (2.4), we choose a positive constant $B$ such that $B > A$, and introduce a new variable $U(\phi)$ through the following *quadratization formula*, that is

$$U(\phi) = \sqrt{F(\phi) + B}. \tag{3.1}$$

Since $F(\phi) + B > -A + B > 0$, we denote

$$H(\phi) = 2\frac{d}{d\phi}U(\phi) = \frac{f(\phi)}{\sqrt{F(\phi) + B}}, \tag{3.2}$$

then the Cahn-Hilliard equation (2.5)-(2.6) can be rewritten as:

$$\phi_t - M\Delta w = 0, \tag{3.3}$$

$$w = -\epsilon^2 \Delta \phi + H(\phi)U, \tag{3.4}$$

$$U_t = \frac{1}{2}H(\phi)\phi_t, \tag{3.5}$$

with the initial conditions

$$\phi|_{t=0} = \phi_0, \ U|_{t=0} = \sqrt{F(\phi_0) + B}. \tag{3.6}$$

Note (3.5) is actually an ODE for the new variable $U$, therefore, the boundary conditions for the new system (3.3)-(3.5) are still (2.7).

The new transformed system (3.3)-(3.5) also follows an energy dissipative law in terms of $\phi$ and the new variable $U$. By taking the $L^2$ inner product of (3.3) with $-w$, of (3.4) with $\phi_t$, of (3.5) with $-2U$, performing integration by parts and summing up all obtained equalities, we can obtain the energy dissipation law of the new system (3.3)-(3.5), that reads as

$$\frac{d}{dt}E(\phi, U) = -M\|\nabla w\|^2, \tag{3.7}$$

where

$$E(\phi, U) = \int_\Omega (\frac{\epsilon^2}{2}|\nabla\phi|^2 + U^2)d\boldsymbol{x}. \tag{3.8}$$

Note the transformed system (3.3)-(3.5) is exactly equivalent to the original system (2.5)-(2.6) since (3.1) can be easily obtained by integrating (3.5) with respect to the time. Therefore, the energy law (3.7) for the transformed system is exactly the same as the energy law (2.10) for the original system for the time-continuous case. We emphasize that we will develop energy stable numerical schemes for time stepping of the new transformed system (3.3)-(3.5). Consequently, the proposed schemes should follow the new energy dissipation law (3.7) formally instead of the energy law for the original system (2.10).

Let $\delta t > 0$ denote the time step size and set $t_n = n\delta t$ for $0 \leq n \leq N$ with $T = N\delta t$, then a first-order, semi-discrete time discretization IEQ scheme for solving the new transformed system (3.3)-(3.5) reads as,

$$\frac{\phi^{n+1} - \phi^n}{\delta t} - M\Delta w^{n+1} = 0, \tag{3.9}$$

$$w^{n+1} = -\epsilon^2 \Delta \phi^{n+1} + H^n U^{n+1}, \tag{3.10}$$

$$U^{n+1} - U^n = \frac{1}{2}H^n(\phi^{n+1} - \phi^n), \tag{3.11}$$

where $H^n = H(\phi^n)$. The boundary conditions are as follows,

$$(i) \ \phi^{n+1}, w^{n+1} \ are \ periodic; \text{or} \ (ii) \ \partial_{\boldsymbol{n}}\phi^{n+1}|_{\partial\Omega} = \partial_{\boldsymbol{n}}w^{n+1}|_{\partial\Omega} = 0. \tag{3.12}$$

We show the unconditionally energy stablilty of the scheme (3.9)-(3.11) as follows.

**Theorem** 3.1. *The scheme (3.9)-(3.11) is unconditionally energy stable in the sense that*

$$E(\phi^{n+1}, U^{n+1}) \leq E(\phi^n, U^n) - \delta t M \|\nabla w^{n+1}\|^2. \tag{3.13}$$

*Proof.* By taking the $L^2$ inner product of (3.9) with $-\delta t w^{n+1}$ and performing integration by parts, we derive

$$-(\phi^{n+1} - \phi^n, w^{n+1}) - M\delta t \|\nabla w^{n+1}\|^2 = 0. \tag{3.14}$$

By taking the $L^2$ inner product of (3.10) with $\phi^{n+1} - \phi^n$, using the identity of

$$2(a, a-b) = a^2 - b^2 + (a-b)^2, \tag{3.15}$$

and performing integration by parts, we get

$$(\phi^{n+1} - \phi^n, w^{n+1}) = \frac{\epsilon^2}{2}(\|\nabla \phi^{n+1}\|^2 - \|\nabla \phi^n\|^2 + \|\nabla \phi^{n+1} - \nabla \phi^n\|^2) + (H^n U^{n+1}, \phi^{n+1} - \phi^n).$$

By taking the $L^2$ inner product of (3.11) with $-2U^{n+1}$ and using (3.15), we get

$$-(\|U^{n+1}\|^2 - \|U^n\|^2 + \|U^{n+1} - U^n\|^2) = -(H^n(\phi^{n+1} - \phi^n), U^{n+1}).$$

By combining the above equations together, we have

$$\frac{\epsilon^2}{2}(\|\nabla \phi^{n+1}\|^2 - \|\nabla \phi^n\|^2 + \|\nabla \phi^{n+1} - \nabla \phi^n\|^2) + \|U^{n+1}\|^2 - \|U^n\|^2 + \|U^{n+1} - U^n\|^2 = -\delta t M \|\nabla w^{n+1}\|^2,$$

which concludes the energy stability (3.13) by dropping some unnecessary positive terms. □

**Remark 3.1.** *The essential idea of the IEQ method is to transform the complicated nonlinear potentials into a simple quadratic form in terms of some new variables via a change of variables. When the nonlinear potential is the double well potential, one can choose $B = 0$ and thus this method is exactly the same as the so-called Lagrange multiplier method developed in [15]. We remark that the Lagrange multiplier method in [15] only works for the fourth order polynomial potential $\phi^4$ since its derivative $\phi^3$ can be decomposed into $\lambda(\phi)\phi$ with $\lambda(\phi) = \phi^2$ which can be viewed as a Lagrange multiplier term. However, for other type potentials, the Lagrange multiplier method is not applicable. For example, one cannot separate a factor of $\phi$ from the logarithmic term.*

*On the contrary, the IEQ approach can handle various complex nonlinear terms as long as the corresponding nonlinear potentials are bounded from below, see the author's recent works in [6, 8, 15, 17, 36, 37, 39–44, 48, 49]). This simple way of quadratization provides some great advantages including (i) the complicated nonlinear potential is transformed into a quadratic polynomial that is much easier to handle; (ii) the derivative of the quadratic polynomial is linear, which makes it possible to develop linear scheme; and (iii) the quadratic formulation in terms of new variables can automatically maintain this property of positivity (or bounded from below) of nonlinear potentials.*

**3.2. Implementations and well-posedness.** The introduction of the new variable $U$ may lead to an increase of the computational cost if one attempts to solve the coupled system (3.9)-(3.11). In fact, note the nonlinear coefficient $H$ of the new variable $U$ is treated explicitly in (3.10), thus we can rewrite it as follows,

$$U^{n+1} = \frac{1}{2}H^n \phi^{n+1} + U^n - \frac{1}{2}H^n \phi^n. \tag{3.16}$$

Using this equality, (3.10)-(3.11) can be rewritten as

$$\phi^{n+1} - \delta t M \Delta w^{n+1} = \phi^n, \tag{3.17}$$

$$-w^{n+1} + P(\phi^{n+1}) = -H^n U^n + \frac{1}{2}H^n H^n \phi^n, \tag{3.18}$$

where $P(\phi) = -\epsilon^2 \Delta \phi + \frac{1}{2} H^n H^n \phi$. Therefore, in practice one can solve $\phi^{n+1}$ and $w^{n+1}$ directly from (3.17)-(3.18). Once $\phi^{n+1}$ is obtained, $U^{n+1}$ is automatically given by (3.16).

The scheme (3.18) includes a variable-coefficient $\frac{1}{2}H^n H^n$ for the term $\phi^{n+1}$ that leads to time-dependent dense matrices. Explicitly building those time-dependent dense matrices are extremely expensive (Note that, if one use finite element methods, the corresponding matrices will be sparse but time-dependent). So in practice, an efficient way is to use a conjugate gradient type solver with preconditioning (PCG), that only needs a subroutine to calculate the matrix-vector product instead of building the mass matrix explicitly. An efficient preconditioner is to replace the variable coefficient term $\frac{1}{2}H^n H^n \phi^{n+1}$ by a term with a constant coefficient being the maximum, i.e., $\frac{1}{2}\|H^n H^n\|_{L^\infty}\phi$.

Furthermore, for any $\phi, \psi$ satisfy the boundary conditions (3.12), we have

$$(3.19) \qquad (P(\phi), \psi) = \frac{\epsilon^2}{2}(\nabla \phi, \nabla \psi) + \frac{1}{2}(H^n \phi, H^n \psi),$$

that means the linear operator $P(\phi)$ is symmetric (self-adjoint). Meanwhile, for $\phi$ that satisfies the boundary conditions (3.12) and $\int_\Omega \phi d\boldsymbol{x} = 0$, we have

$$(3.20) \qquad (P(\phi), \phi) = \frac{\epsilon^2}{2}\|\nabla \phi\|^2 + \frac{1}{2}\|H^n \phi\|^2 \geq 0,$$

where " = " is valid if and only if $\phi \equiv 0$. These facts imply that the linear operator $P$ is symmetric positive definite in $H_{per}(\Omega)$ or $H(\Omega)$.

We now show the well-posedness of the scheme (3.17)-(3.18). In the following arguments, we will only consider the periodic boundary condition for convenience. For the case of homogenous Neumann boundary conditions, as long as the space of $H_{per}(\Omega)$ is replaced by $H(\Omega)$, all theoretical results are still valid.

By taking the $L^2$ inner product of (3.17) with 1, we derive

$$(3.21) \qquad \int_\Omega \phi^{n+1} d\boldsymbol{x} = \int_\Omega \phi^n d\boldsymbol{x} = \cdots = \int_\Omega \phi^0 d\boldsymbol{x}.$$

Let $V_\phi = \frac{1}{|\Omega|}\int_\Omega \phi^0 d\boldsymbol{x}$, $V_w = \frac{1}{\Omega}\int_\Omega \mu^{n+1}$, and we define

$$(3.22) \qquad \widehat{\phi}^{n+1} = \phi^{n+1} - V_\phi, \widehat{w}^{n+1} = w^{n+1} - V_w,$$

then the weak form for (3.17)-(3.18) is the system with unknowns $\phi, w \in (H_{per}, H_{per})(\Omega)$:

$$(3.23) \qquad (\phi, \mu) + \delta t M(\nabla w, \nabla \mu) = (\widehat{\phi}^n, \mu), \quad \mu \in H_{per}(\Omega),$$

$$(3.24) \quad (-w, \psi) + \epsilon^2(\nabla \phi, \nabla \psi) + \frac{1}{2}(H^n \phi, H^n \psi) = (-H^n U^n + \frac{1}{2}H^n H^n \widehat{\phi}^n, \psi), \quad \psi \in H_{per}(\Omega).$$

We denote the above linear system as

$$(3.25) \qquad (\boldsymbol{L}\boldsymbol{X}, \boldsymbol{Y}) = (\boldsymbol{B}, \boldsymbol{Y}),$$

where $\boldsymbol{X} = (w, \phi)^T$, $\boldsymbol{Y} = (\mu, \psi)^T$, $\boldsymbol{B} = (\widehat{\phi}^n, -H^n U^n + \frac{1}{2}H^n H^n \widehat{\phi}^n)^T$, and $\boldsymbol{X}, \boldsymbol{Y} \in (H_{per}, H_{per})(\Omega)$.

The well-posedness of the linear system (3.25) is shown as follows.

**Theorem** 3.2. *The linear system* (3.25) *admits a unique solution* $(w, \phi)$ *in* $(H_{per}, H_{per})(\Omega)$.

*Proof.* (i) For any $\boldsymbol{X} = (w, \phi)^T$, $\boldsymbol{Y} = (\mu, \psi)^T$ with $\boldsymbol{X}, \boldsymbol{Y} \in (H_{per}, H_{per})(\Omega)$, we have

$$(3.26) \qquad (\boldsymbol{L}\boldsymbol{X}, \boldsymbol{Y}) \leq c_1(\|\phi\|_1 + \|w\|_1)(\|\psi\|_1 + \|\mu\|_1),$$

where $c_1$ is a postive constant dependent on $\delta t, M, \epsilon$ and $\|H^n\|_{L^\infty}$. Therefore, the bilinear form $(\boldsymbol{L}\boldsymbol{X}, \boldsymbol{Y})$ is bounded.

(ii) For any $\boldsymbol{X} = (w, \phi)^T \in (H_{per}, H_{per})(\Omega)$, we derive

$$(\boldsymbol{LX}, \boldsymbol{X}) = \delta t M \|\nabla w\|^2 + \epsilon^2 \|\nabla \phi\|^2 + \frac{1}{2}\|H^n \phi\|^2 \geq c_2(\|w\|_1^2 + \|\phi\|_1^2), \tag{3.27}$$

where $c_2$ is a constant dependent on $\delta t, M, \epsilon$. Thus the bilinear form $(\boldsymbol{LX}, \boldsymbol{Y})$ is coercive. Then from the Lax-Milgram Theorem, we conclude that the linear system (3.25) admits a unique solution $(w, \phi) \in (H_{per}, H_{per})(\Omega)$. Namely, the scheme (3.17)-(3.18) admits a unique solution $(w^{n+1}, \phi^{n+1}) \in (H^1, H^1)(\Omega)$. □

**3.3. Error estimates.** We now focus on the error estimates. To simplify the notations, without loss of generality, in the below, we let $M = \epsilon = 1$. We use $x \lesssim y$ to denote there exists a constant $C$ that is independent on $\delta t$ and $n$ such that $x \leq Cy$.

To derive the error estimates of (3.9)-(3.11), we first give the following two Lemmas to establish the quantitative relation between the $L^2$ and $H^1$ norms of $H(\phi(t^n)) - H(\phi^n)$ and $\phi(t^n) - \phi^n$ under some reasonable assumptions.

**Lemma** 3.1. *Suppose (i) $F(x)$ is uniformly bounded from below: $F(x) > -A$ for any $x \in (-\infty, \infty)$; (ii) $F(x) \in C^2(-\infty, \infty)$; and (iii) there exists a positive constant $C_0$ such that*

$$\max_{n \leq M}(\|\phi(t_n)\|_{L^\infty}, \|\phi^n\|_{L^\infty}) \leq C_0, \tag{3.28}$$

*then we have*

$$\|H(\phi(t_n)) - H^n\| \leq \widehat{C}_0 \|\phi(t_n) - \phi^n\|, \tag{3.29}$$

*for $n \leq M$, where $\widehat{C}_0$ is a positive constant that is only dependent of $C_0, A$ and $B$.*

*Proof.* First, for any $\theta \in [0, 1]$, (3.28) can ensure $\psi^n = \theta \phi(t_n) + (1-\theta)\phi^n$ is uniformly bounded, i.e., $\psi^n \in [-2C_0, 2C_0]$ for $n \leq M$. Thus from assumption (ii), we can always find a positive constant $C_1$ such that

$$\max_{n \leq M}\left(\|F(\psi^n)\|_{L^\infty}, \|f(\psi^n)\|_{L^\infty}, \|f'(\psi^n)\|_{L^\infty}, \|\sqrt{F(\psi^n)+B}\|_{L^\infty}\right) \leq C_1. \tag{3.30}$$

Second, for any $x, y \in (-A, \infty)$, by applying the intermediate value Theorem, there exists some value $\xi \in (-A, \infty)$ that is between $x$ and $y$, such that $\sqrt{x+B} - \sqrt{y+B} = \frac{1}{2\sqrt{\xi+B}}(x-y)$, that implies

$$|\sqrt{x+B} - \sqrt{y+B}| \leq \frac{1}{2\sqrt{B-A}}|x-y|. \tag{3.31}$$

Thus, using (3.30), (3.31), and applying the intermediate value Theorem again, we derive

$$\begin{aligned}\left|\sqrt{F(\phi(t_n))+B} - \sqrt{F(\phi^n)+B}\right| &\leq \frac{1}{2\sqrt{B-A}}|F(\phi(t_n)) - F(\phi^n)| \\ &= \frac{1}{2\sqrt{B-A}}|f(\theta\phi(t^n)+(1-\theta)\phi^n)||\phi(t^n)-\phi^n| \\ &\leq \frac{C_1}{2\sqrt{B-A}}|\phi(t^n)-\phi^n|.\end{aligned} \tag{3.32}$$

Third, for $n \leq M$, from (3.30) and (3.32), we derive

$$
\begin{aligned}
|H(\phi(t_n)) - H^n| &= \left| \frac{f(\phi(t_n))}{\sqrt{F(\phi(t_n)) + B}} - \frac{f(\phi^n)}{\sqrt{F(\phi^n) + B}} \right| \\
&= \frac{\left| f(\phi(t_n))\sqrt{F(\phi^n) + B} - f(\phi^n)\sqrt{F(\phi(t_n)) + B} \right|}{\sqrt{F(\phi(t_n)) + B}\sqrt{F(\phi^n) + B}} \\
&\leq \frac{1}{B - A} \left| f(\phi(t_n))\sqrt{F(\phi^n) + B} - f(\phi^n)\sqrt{F(\phi(t_n)) + B} \right| \\
&\leq \frac{1}{B - A} |f(\phi(t_n))| \left| \sqrt{F(\phi(t_n)) + B} - \sqrt{F(\phi^n) + B} \right| \\
&\quad + \frac{1}{B - A} \sqrt{F(\phi(t_n)) + B} \, |f(\phi^n) - f(\phi(t_n))| \\
&\leq \frac{C_1}{B - A} \frac{C_1}{2\sqrt{B - A}} |\phi(t^n) - \phi^n| + \frac{C_1}{B - A} |f'(\psi^n)||\phi(t_n) - \phi^n| \\
&\leq \left( \frac{C_1^2}{2\sqrt{(B - A)^3}} + \frac{C_1^2}{B - A} \right) |\phi(t^n) - \phi^n|,
\end{aligned}
$$
(3.33)

where we have used $\frac{1}{\sqrt{F(x) + B}} \leq \frac{1}{\sqrt{B - A}}$ for any $x \in (-\infty, \infty)$.

Finally, let $\widehat{C}_0 = \frac{C_1^2}{2\sqrt{(B-A)^3}} + \frac{C_1^2}{B-A}$, we derive

$$\|H(\phi(t_n)) - H^n\| = \left( \int_\Omega |H(\phi(t_n)) - H(\phi^n)|^2 d\boldsymbol{x} \right)^{\frac{1}{2}} \leq \widehat{C}_0 \left( \int_\Omega |\phi(t_n) - \phi^n|^2 d\boldsymbol{x} \right)^{\frac{1}{2}} = \widehat{C}_0 \|\phi(t_n) - \phi^n\|.$$

**Remark 3.2.** *Lemma 3.1 establishes a quantitative relation between the $L^2$ norm of $H(\phi(t^n)) - H(\phi^n)$ and $\phi(t^n) - \phi^n$ under some reasonable assumptions, where the Lipschitz property (3.31) of the quadratization function $\sqrt{x + B}$ ($x > -A$) plays a critical role. We note assumptions (i) and (ii) are automatically valid for the fourth order polynomial type double-well potential. Indeed, for the double well potential, one can choose $B = 0$ and $H^n = \phi^n$, thus Lemma 3.1 is trivial and the error analysis is straightforward for this case, see [7, 15, 16]. However, for the logarithmic Flory-Huggins potential, $B \neq 0$ and assumption (ii) is not true since the domain is the open interval $(0, 1)$ instead of $(-\infty, \infty)$. This issue can be overcome by extending the logarithmic functional near the domain boundary with a continuous, convex, piecewise function, see [9, 11, 36]. Such a regularized method is also a common practice to remove the difficulty about that any small fluctuation near the domain boundary $(0, 1)$ can cause the overflow, numerically. In this way, the domain is regularized to $(-\infty, \infty)$ and thus the assumptions (ii) will become valid.*

Similarly, we further establish the relation between their $H^1$ norms, as follows.

**Lemma** 3.2. *Suppose (i) $F(x)$ is uniformly bounded from below: $F(x) > -A$ for any $x \in (-\infty, \infty)$; (ii) $F(x) \in C^3(-\infty, \infty)$; and (iii) there exists a positive constant $D_0$ such that*

$$\max_{n \leq M}(\|\phi(t_n)\|_{L^\infty}, \|\phi^n\|_{L^\infty}, \|\nabla\phi(t_n)\|_{L^\infty}) \leq D_0, \tag{3.34}$$

*then we have*

$$\|\nabla H(\phi(t_n)) - \nabla H^n\| \leq \widehat{D}_0(\|\phi(t_n) - \phi^n\| + \|\nabla\phi(t_n) - \nabla\phi^n\|), \tag{3.35}$$

*for $n \leq M$, where $\widehat{D}_0$ is a positive constant dependent on $D_0, A$ and $B$.*

*Proof.* First, from assumption (ii) and (iii), for any $\psi^n = \theta\phi(t_n) + (1 - \theta)\phi^n$ where $\theta \in [0, 1]$, we can always find a positive constant $D_1$ that is dependent on $D_0$, such that

$$(3.36) \quad \max_{n \leq M} \left( \|F(\psi^n)\|_{L^\infty}, \|f(\psi^n)\|_{L^\infty}, \|f'(\psi^n)\|_{L^\infty}, \|f''(\psi^n)\|_{L^\infty}, \|F(\psi^n) + B\|_{L^\infty} \right) \leq D_1.$$

Second, for convenience, we denote $G(u) = f'(u)(F(u) + B) - \frac{1}{2}f^2(u)$. From (3.36), assumption (i), (ii) and (iii), we derive

$$(3.37) \quad \begin{aligned} |\nabla H(\phi(t_n)) - \nabla H(\phi^n)| &= |H'(\phi(t_n))\nabla\phi(t_n) - H'(\phi^n)\nabla\phi^n| \\ &\leq |\nabla\phi(t_n)||H'(\phi(t_n)) - H'(\phi^n)| + |H'(\phi^n)||\nabla\phi(t_n) - \nabla\phi^n| \\ &\leq D_0|H'(\phi(t_n)) - H'(\phi^n)| + |H'(\phi^n)||\nabla\phi(t_n) - \nabla\phi^n| \\ &\leq D_0|H'(\phi(t_n)) - H'(\phi^n)| + \left|\frac{G(\phi^n)}{(F(\phi^n) + B)^{\frac{3}{2}}}\right| |\nabla\phi(t_n) - \nabla\phi^n| \\ &\leq D_0|H'(\phi(t_n)) - H'(\phi^n)| + \frac{\frac{3}{2}D_1^2}{(B - A)^{\frac{3}{2}}}|\nabla\phi(t_n) - \nabla\phi^n|, \end{aligned}$$

where we have used $|G(\phi^n)| \leq \frac{3}{2}D_1^2$ from (3.36) and $\frac{1}{(F(\phi^n)+B)^{\frac{3}{2}}} \leq \frac{1}{(B-A)^{\frac{3}{2}}}$ from assumption (i).

Third, we estimate

$$(3.38) \quad \begin{aligned} |H'(\phi(t_n)) - H'(\phi^n)| &= \left| \frac{G(\phi(t_n))}{(F(\phi(t_n)) + B)^{\frac{3}{2}}} - \frac{G(\phi^n)}{(F(\phi^n) + B)^{\frac{3}{2}}} \right| \\ &= \frac{\left|(F(\phi^n) + B)^{\frac{3}{2}}G(\phi(t_n)) - (F(\phi(t_n)) + B)^{\frac{3}{2}}G(\phi^n)\right|}{(F(\phi(t_n)) + B)^{\frac{3}{2}}(F(\phi^n) + B)^{\frac{3}{2}}} \\ &\leq \frac{1}{(B - A)^3} \left|(F(\phi^n) + B)^{\frac{3}{2}}G(\phi(t_n)) - (F(\phi(t_n)) + B)^{\frac{3}{2}}G(\phi^n)\right| \\ &\leq \frac{1}{(B - A)^3} \left|(F(\phi^n) + B)^{\frac{3}{2}} - (F(\phi(t_n)) + B)^{\frac{3}{2}}\right| |G(\phi(t_n))| \\ &\quad + \frac{1}{(B - A)^3}(F(\phi(t_n)) + B)^{\frac{3}{2}} |f'(\phi(t_n))(F(\phi(t_n)) + B) - f'(\phi^n)(F(\phi^n) + B)| \\ &\quad + \frac{1}{(B - A)^3}(F(\phi(t_n)) + B)^{\frac{3}{2}}\frac{1}{2}\left|f^2(\phi^n) - f^2(\phi(t_n))\right| \\ &\leq \frac{\frac{3}{2}D_1^2}{(B - A)^3}\left|(F(\phi^n) + B)^{\frac{3}{2}} - (F(\phi(t_n)) + B)^{\frac{3}{2}}\right| \qquad (: \text{term } I_1) \\ &\quad + \frac{D_1^{\frac{3}{2}}}{(B - A)^3} |f'(\phi(t_n))(F(\phi(t_n)) + B) - f'(\phi^n)(F(\phi^n) + B)| \qquad (: \text{term } I_2) \\ &\quad + \frac{1}{2}\frac{D_1^{\frac{3}{2}}}{(B - A)^3} \left|f^2(\phi^n) - f^2(\phi(t_n))\right| \qquad (: \text{term } I_3). \end{aligned}$$

For term $I_1$, by applying the intermediate value Theorem twice and using (3.36), we derive

$$(3.39) \quad \begin{aligned} I_1 &\leq \frac{\frac{3}{2}D_1^2}{(B - A)^3}\frac{3}{2}\sqrt{(\xi_1 + B)}\left|F(\phi^n) - F(\phi(t^n))\right| \\ &\leq \frac{\frac{3}{2}D_1^2}{(B - A)^3}\frac{3}{2}\sqrt{(\xi_1 + B)}|f(\xi_2)|\left|\phi^n - \phi(t^n)\right| \\ &\leq \frac{\frac{9}{4}D_1^3\sqrt{(2D_1 + B)}}{(B - A)^3}\left|\phi^n - \phi(t^n)\right|, \end{aligned}$$

where $\xi_1 = \theta_1 F(\phi(t^n)) + (1 - \theta_1)F(\phi^n)$, $\xi_2 = \theta_2 \phi(t^n) + (1 - \theta_2)\phi^n$ for some $\theta_1, \theta_2 \in [0, 1]$, $\sqrt{\xi_1 + B} \leq$

$\sqrt{2D_1 + B}$ and $f(\xi_2) \leq D_1$.

For term $I_2$, using the intermediate value Theorem for $F$ and $f'$ and (3.36), we derive

$$
\begin{aligned}
I_2 &\leq \frac{D_1^{\frac{3}{2}}}{(B-A)^3}\Big(|f'(\phi(t_n))(F(\phi(t_n)) - F(\phi^n))| + |(f'(\phi(t_n)) - f'(\phi^n))(F(\phi^n) + B)|\Big) \\
&\leq \frac{D_1^{\frac{5}{2}}}{(B-A)^3}\Big(|(F(\phi(t_n)) - F(\phi^n))| + |f'(\phi(t_n)) - f'(\phi^n)|\Big) \\
&\leq \frac{D_1^{\frac{7}{2}}}{(B-A)^3}|\phi^n - \phi(t^n)|.
\end{aligned}
\tag{3.40}
$$

For term $I_3$, using the intermediate value Theorem for $f$ and (3.36), we derive

$$
\begin{aligned}
I_3 &\leq \frac{1}{2}\frac{D_1^{\frac{3}{2}}}{(B-A)^3}|f(\phi(t^n)) + f(\phi^n)||f(\phi(t^n)) - f(\phi^n)| \\
&\leq \frac{D_1^{\frac{5}{2}}}{(B-A)^3}|f(\phi(t^n)) - f(\phi^n)| \\
&\leq \frac{D_1^{\frac{7}{2}}}{(B-A)^3}|\phi^n - \phi(t^n)|.
\end{aligned}
\tag{3.41}
$$

Thus, by combining (3.39), (3.40), (3.41), and denote $\widehat{D}_1 = \frac{\frac{9}{4}D_1^3\sqrt{(2D_1+B)}}{(B-A)^3} + \frac{2D_1^{\frac{7}{2}}}{(B-A)^3}$, we have

$$
|H'(\phi(t^n)) - H'(\phi^n)| \leq \widehat{D}_1|\phi^n - \phi(t^n)|. \tag{3.42}
$$

Therefore, from (3.37) and let $D_2 = max(D_0\widehat{D}_1, \frac{\frac{3}{2}D_1^2}{(B-A)^{\frac{3}{2}}})$, we derive

$$
|\nabla H(\phi(t_n)) - \nabla H(\phi^n)| \leq D_2\Big(|\phi(t_n) - \phi^n| + |\nabla\phi(t_n) - \nabla\phi^n|\Big).
$$

Then we have

$$
\begin{aligned}
\|\nabla H(\phi(t_n)) - \nabla H(\phi^n)\|^2 &= \int_\Omega |\nabla H(\phi(t_n)) - \nabla H(\phi^n)|^2 dx \\
&\leq 2D_2^2 \int_\Omega |\phi(t_n) - \phi^n|^2 + |\nabla\phi(t_n) - \nabla\phi^n|^2 dx \\
&= \widehat{D}_0^2\Big(\|\phi(t_n) - \phi^n\|^2 + \|\nabla\phi(t_n) - \nabla\phi^n\|^2\Big),
\end{aligned}
\tag{3.43}
$$

where $\widehat{D}_0 = \sqrt{2D_2}$, that concludes (3.35). □

We now establish the error estimates for scheme (3.9)-(3.11). To this end, we formulate the Cahn-Hilliard system (3.3)-(3.5) as a truncation form:

$$
\frac{\phi(t_{n+1}) - \phi(t_n)}{\delta t} = \Delta w(t_{n+1}) + R_\phi^{n+1}, \tag{3.44}
$$

$$
w(t_{n+1}) = -\Delta\phi(t_{n+1}) + H(\phi(t_n))U(t_{n+1}) + R_w^{n+1}, \tag{3.45}
$$

$$
U(t_{n+1}) - U(t_n) = \frac{1}{2}H(\phi(t_n))(\phi(t_{n+1}) - \phi(t_n)) + \delta t R_u^{n+1}, \tag{3.46}
$$

where

$$(3.47) \begin{cases} R_\phi^{n+1} = \dfrac{\phi(t_{n+1}) - \phi(t_n)}{\delta t} - \phi_t(t_{n+1}), \\ R_w^{n+1} = H(\phi(t_{n+1}))U(t_{n+1}) - H(\phi(t_n))U(t_{n+1}), \\ R_u^{n+1} = \dfrac{U(t_{n+1}) - U(t_n)}{\delta t} - U_t(t_{n+1}) + \dfrac{1}{2}H(\phi(t_{n+1}))\phi_t(t_{n+1}) - \dfrac{1}{2}H(\phi(t_n))\dfrac{\phi(t_{n+1}) - \phi(t_n)}{\delta t}. \end{cases}$$

We assume the exact solution $\phi, w, U$ of the system (3.3)-(3.5) possesses the following regularity conditions,

$$(3.48) \begin{cases} \phi \in L^\infty(0,T;H^2(\Omega)) \cap L^\infty(0,T;W^{1,\infty}(\Omega)), \\ U \in L^\infty(0,T;W^{1,\infty}(\Omega)), \\ w \in L^\infty(0,T;H^1(\Omega)), \\ \phi_t \in L^2(0,T;H^1(\Omega)) \cap L^\infty(0,T;L^\infty), \phi_{tt}, U_{tt} \in L^2(0,T;L^2(\Omega)). \end{cases}$$

One can easily establish the following estimates for the truncation errors, provided that the exact solutions of the system (3.3)-(3.5) satisfy (3.48).

**Lemma 3.3.** *Under the regularity conditions* (3.48), *the truncation errors satisfy*

$$(3.49) \qquad \delta t \sum_{n=0}^{[\frac{T}{\delta t}]} (\|R_\phi^{n+1}\|^2 + \|R_w^{n+1}\|_1^2 + \|R_u^{n+1}\|^2) \lesssim \delta t^2.$$

*Proof.* Since the proof is rather straight forward, we leave this to the interested readers. □

To derive the error estimates, we denote the error functions as

$$(3.50) \begin{cases} e_\phi^n = \phi(t_n) - \phi^n, & e_w^n = w(t_n) - w^n, \\ e_u^n = U(t_n) - U^n, & e_H^n = H(\phi(t_n)) - H(\phi^n). \end{cases}$$

By subtracting (3.9)-(3.11) from (3.44)-(3.46), we derive the error equations:

$$(3.51) \qquad \frac{e_\phi^{n+1} - e_\phi^n}{\delta t} = \Delta e_w^{n+1} + R_\phi^{n+1},$$

$$(3.52) \qquad e_w^{n+1} = -\Delta e_\phi^{n+1} + e_H^n U(t_{n+1}) + H^n e_u^{n+1} + R_w^{n+1},$$

$$(3.53) \qquad e_u^{n+1} - e_u^n = \frac{1}{2}(e_H^n(\phi(t_{n+1}) - \phi(t_n)) + H^n(e_\phi^{n+1} - e_\phi^n)) + \delta t R_u^{n+1}.$$

We first prove the $L^\infty$ stability of $\phi^n$, which plays the key point in the error estimates. From (3.48), we define a positive constant $\kappa$ such that

$$(3.54) \qquad \kappa = \max_{0 \le t \le T} \|\phi(t)\|_{L^\infty} + 1.$$

The preliminary result is given in the following lemma.

**Lemma 3.4.** *Suppose (i) $F(x)$ is uniformly bounded from below: $F(x) > -A$ for any $x \in (-\infty, \infty)$; (ii) $F(x) \in C^3(-\infty, \infty)$; and (iii) the exact solutions of (3.3)-(3.5) satisfy the regularity conditions (3.48), then there exists a positive constant $s_0$ that is given in the proof, such that when $\delta t \le s_0$, the solution $\phi^n$ of (3.9)-(3.11) is uniformly bounded as*

$$(3.55) \qquad \|\phi^n\|_{L^\infty} \le \kappa, n = 0, 1, \cdots, N = \frac{T}{\delta t}.$$

*Proof.* We use the mathematical induction to prove this Lemma.

For $n = 0$, $\|\phi^0\|_{L^\infty} \leq \kappa$ is true naturally.

Assuming that $\|\phi^n\|_{L^\infty} \leq \kappa$ is valid for all $n \leq M$, we show $\|\phi^{M+1}\|_{L^\infty} \leq \kappa$ is also valid through the following two steps.

(i). By taking the $L^2$ inner product of (3.51) with $\delta t e_w^{n+1}$, we obtain

$$(e_\phi^{n+1} - e_\phi^n, e_w^{n+1}) + \delta t \|\nabla e_w^{n+1}\|^2 = \delta t (R_\phi^{n+1}, e_w^{n+1}). \tag{3.56}$$

By taking the $L^2$ inner product of (3.52) with $-(e_\phi^{n+1} - e_\phi^n)$, we obtain

$$\begin{aligned}-(e_w^{n+1}, e_\phi^{n+1} - e_\phi^n) &+ \frac{1}{2}(\|\nabla e_\phi^{n+1}\|^2 - \|\nabla e_\phi^n\|^2 + \|\nabla e_\phi^{n+1} - \nabla e_\phi^n\|^2) \\ &= -(e_H^n U(t_{n+1}) + H^n e_u^{n+1}, e_\phi^{n+1} - e_\phi^n) - (R_w^{n+1}, e_\phi^{n+1} - e_\phi^n).\end{aligned} \tag{3.57}$$

By taking the $L^2$ inner product of (3.53) with $2e_u^{n+1}$, we get

$$\begin{aligned}\|e_u^{n+1}\|^2 &- \|e_u^n\|^2 + \|e_u^{n+1} - e_u^n\|^2 \\ &= (e_H^n(\phi(t_{n+1}) - \phi(t_n)) + H^n(e_\phi^{n+1} - e_\phi^n), e_u^{n+1}) + 2\delta t(R_u^{n+1}, e_u^{n+1}).\end{aligned} \tag{3.58}$$

By taking the $L^2$ inner product of (3.51) with $\delta t e_\phi^{n+1}$, we get

$$\frac{1}{2}(\|e_\phi^{n+1}\|^2 - \|e_\phi^n\|^2 + \|e_\phi^{n+1} - e_\phi^n\|^2) = \delta t(\Delta e_w^{n+1}, e_\phi^{n+1}) + \delta t(R_\phi^{n+1}, e_\phi^{n+1}). \tag{3.59}$$

By taking the $L^2$ inner product of (3.52) with $\delta t e_w^{n+1}$, we derive

$$\delta t \|e_w^{n+1}\|^2 = -\delta t(\Delta e_\phi^{n+1}, e_w^{n+1}) + \delta t(e_H^n U(t_{n+1}) + H^n e_u^{n+1}, e_w^{n+1}) + \delta t(R_w^{n+1}, e_w^{n+1}). \tag{3.60}$$

Combining (3.56)-(3.60) together, we obtain

$$\begin{aligned}&\frac{1}{2}(\|e_\phi^{n+1}\|^2 + \|\nabla e_\phi^{n+1}\|^2 - \|e_\phi^n\|^2 - \|\nabla e_\phi^n\|^2) + (\|e_u^{n+1}\|^2 - \|e_u^n\|^2) + \delta t(\|e_w^{n+1}\|^2 + \|\nabla e_w^{n+1}\|^2) \\ &+ \frac{1}{2}(\|e_\phi^{n+1} - e_\phi^n\|^2 + \|\nabla e_\phi^{n+1} - \nabla e_\phi^n\|^2) + \|e_u^{n+1} - e_u^n\|^2 \\ &= -(e_H^n U(t_{n+1}), e_\phi^{n+1} - e_\phi^n) + (e_H^n(\phi(t_{n+1}) - \phi(t_n)), e_u^{n+1}) + \delta t(e_H^n U(t_{n+1}) + H^n e_u^{n+1}, e_w^{n+1}) \\ &- (R_w^{n+1}, e_\phi^{n+1} - e_\phi^n) + \delta t(R_\phi^{n+1}, e_w^{n+1}) + 2\delta t(R_u^{n+1}, e_u^{n+1}) + \delta t(R_\phi^{n+1}, e_\phi^{n+1}) + \delta t(R_w^{n+1}, e_w^{n+1}).\end{aligned} \tag{3.61}$$

By using Lemma 3.1, Lemma 3.2 and (3.48), for $n \leq M$, we estimate terms on the right hand side as follows.

$$\begin{aligned}|(e_H^n U(t_{n+1}), e_\phi^{n+1} - e_\phi^n)| &= \delta t |(e_H^n U(t_{n+1}), \frac{e_\phi^{n+1} - e_\phi^n}{\delta t})| \\ &= \delta t |(e_H^n U(t_{n+1}), \Delta e_w^{n+1} + R_\phi^{n+1})| \\ &= \delta t |(\nabla(e_H^n U(t_{n+1})), \nabla e_w^{n+1}) + (e_H^n U(t_{n+1}), R_\phi^{n+1})| \\ &= \delta t |(U(t_{n+1})\nabla e_H^n + e_H^n \nabla U(t_{n+1}), \nabla e_w^{n+1}) + (e_H^n U(t_{n+1}), R_\phi^{n+1})| \\ &\leq \delta t \|U(t_{n+1})\|_{L^\infty}\|\nabla e_H^n\|\|\nabla e_w^{n+1}\| + \delta t \|e_H^n\|\|\nabla U(t_{n+1})\|_{L^\infty}\|\nabla e_w^{n+1}\| \\ &\quad + \delta t \|e_H^n\|\|U(t_{n+1})\|_{L^\infty}\|R_\phi^{n+1}\| \\ &\lesssim \delta t \|\nabla e_H^n\|\|\nabla e_w^{n+1}\| + \delta t \|e_H^n\|\|\nabla e_w^{n+1}\| + \delta t \|e_H^n\|\|R_\phi^{n+1}\| \\ &\lesssim \frac{1}{4}\delta t\|\nabla e_w^{n+1}\|^2 + \delta t\|\nabla e_H^n\|^2 + \delta t\|e_H^n\|^2 + \delta t\|R_\phi^{n+1}\|^2 \\ &\lesssim \frac{1}{4}\delta t\|\nabla e_w^{n+1}\|^2 + \delta t\|\nabla e_\phi^n\|^2 + \delta t\|e_\phi^n\|^2 + \delta t\|R_\phi^{n+1}\|^2;\end{aligned} \tag{3.62}$$

$$
\begin{aligned}
\left|(e_H^n(\phi(t_{n+1}) - \phi(t_n)), e_u^{n+1})\right| &\leq \|e_H^n\|_{L^4}\|\phi(t_{n+1}) - \phi(t_n)\|_{L^4}\|e_u^{n+1}\| \\
&\lesssim \delta t\|e_H^n\|_{L^4}\|e_u^{n+1}\| \\
&\lesssim \delta t\|e_H^n\|_{L^4}^2 + \delta t\|e_u^{n+1}\|^2 \\
&\lesssim \delta t(\|e_H^n\|^2 + \|\nabla e_H^n\|^2) + \delta t\|e_u^{n+1}\|^2 \\
&\lesssim \delta t(\|e_\phi^n\|^2 + \|\nabla e_\phi^n\|^2) + \delta t\|e_u^{n+1}\|^2;
\end{aligned}
\tag{3.63}
$$

$$
\begin{aligned}
\delta t|(e_H^n U(t_{n+1}) + H^n e_u^{n+1}, e_w^{n+1})| &\leq \delta t(\|e_H^n\|\|U(t_{n+1})\|_{L^\infty} + \|H^n\|_{L^\infty}\|e_u^{n+1}\|)\|e_w^{n+1}\| \\
&\lesssim \delta t(\|e_H^n\| + \|e_u^{n+1}\|)\|e_w^{n+1}\| \\
&\lesssim \frac{1}{6}\delta t\|e_w^{n+1}\|^2 + \delta t\|e_H^n\|^2 + \delta t\|e_u^{n+1}\|^2 \\
&\lesssim \frac{1}{6}\delta t\|e_w^{n+1}\|^2 + \delta t\|e_\phi^n\|^2 + \delta t\|e_u^{n+1}\|^2,
\end{aligned}
\tag{3.64}
$$

where, $\|H^n\|_{L^\infty}$ is bounded since $\|\phi^n\|_{L^\infty}$ is bounded, $f$ is continuous, and $F(x) > -A$;

$$
\begin{aligned}
|(R_w^{n+1}, e_\phi^{n+1} - e_\phi^n)| &= \delta t|(R_w^{n+1}, \frac{e_\phi^{n+1} - e_\phi^n}{\delta t})| \\
&= \delta t|(R_w^{n+1}, \Delta e_w^{n+1} + R_\phi^{n+1})| \\
&= \delta t|(\nabla R_w^{n+1}, \nabla e_w^{n+1}) + (R_w^{n+1}, R_\phi^{n+1})| \\
&\leq \delta t\|\nabla R_w^{n+1}\|\|\nabla e_w^{n+1}\| + \delta t\|R_w^{n+1}\|\|R_\phi^{n+1}\| \\
&\lesssim \frac{1}{4}\delta t\|\nabla e_w^{n+1}\|^2 + \delta t\|\nabla R_w^{n+1}\|^2 + \delta t\|R_w^{n+1}\|^2 + \delta t\|R_\phi^{n+1}\|^2;
\end{aligned}
\tag{3.65}
$$

$$|\delta t(R_\phi^{n+1}, e_w^{n+1})| \leq \delta t\|R_\phi^{n+1}\|\|e_w^{n+1}\| \lesssim \frac{1}{6}\delta t\|e_w^{n+1}\|^2 + \delta t\|R_\phi^{n+1}\|^2; \tag{3.66}$$

$$2\delta t|(R_u^{n+1}, e_u^{n+1})| \lesssim \delta t\|R_u^{n+1}\|^2 + \delta t\|e_u^{n+1}\|^2; \tag{3.67}$$

$$|\delta t(R_\phi^{n+1}, e_\phi^{n+1})| \lesssim \delta t\|R_\phi^{n+1}\|^2 + \delta t\|e_\phi^{n+1}\|^2; \tag{3.68}$$

$$|\delta t(R_w^{n+1}, e_w^{n+1})| \leq \delta t\|R_w^{n+1}\|\|e_w^{n+1}\| \lesssim \frac{1}{6}\delta t\|e_w^{n+1}\|^2 + \delta t\|R_w^{n+1}\|^2. \tag{3.69}$$

Combining the above estimates with (3.61), we obtain

$$
\begin{aligned}
&\|e_\phi^{n+1}\|^2 + \|\nabla e_\phi^{n+1}\|^2 - \|e_\phi^n\|^2 - \|\nabla e_\phi^n\|^2 + 2(\|e_u^{n+1}\|^2 - \|e_u^n\|^2) \\
&\quad + \delta t(\|e_w^{n+1}\|^2 + \|\nabla e_w^{n+1}\|^2) \\
&\lesssim \delta t(\|e_\phi^n\|^2 + \|\nabla e_\phi^n\|^2 + \|e_u^{n+1}\|^2 + \|e_\phi^{n+1}\|^2) \\
&\quad + \delta t(\|R_w^{n+1}\|^2 + \|\nabla R_w^{n+1}\|^2 + \|R_\phi^{n+1}\|^2 + \|R_u^{n+1}\|^2).
\end{aligned}
\tag{3.70}
$$

Summing up the above inequality from $n = 0$ to $m$ $(m \leq M)$ and using Lemma 3.3, we have

$$\|e_\phi^{m+1}\|^2 + \|\nabla e_\phi^{m+1}\|^2 + 2\|e_u^{m+1}\|^2 + \delta t \sum_{n=0}^m (\|e_w^{n+1}\|^2 + \|\nabla e_w^{n+1}\|^2)$$
$$\lesssim \delta t \sum_{n=0}^m (\|e_\phi^n\|^2 + \|\nabla e_\phi^n\|^2 + \|e_u^{n+1}\|^2 + \|e_\phi^{n+1}\|^2) + \delta t \sum_{n=0}^m (\|R_w^{n+1}\|_1^2 + \|R_\phi^{n+1}\|^2 + \|R_u^{n+1}\|^2)$$
$$\lesssim \delta t \sum_{n=0}^m (\|e_\phi^{n+1}\|^2 + \|\nabla e_\phi^{n+1}\|^2 + \|e_u^{n+1}\|^2) + \delta t^2.$$

Then, by using the Gronwall's inequality, there exist two positive constants $s_1, s_2$ such that when $\delta t \leq s_1$, the following inequality holds for any $m \leq M$,

$$(3.71) \qquad \|e_\phi^{m+1}\|^2 + \|\nabla e_\phi^{m+1}\|^2 + \|e_u^{m+1}\|^2 + \delta t \sum_{n=0}^m (\|e_w^{n+1}\|^2 + \|\nabla e_w^{n+1}\|^2) \leq s_2 \delta t^2.$$

(ii). By using the $H^2$ regularity of elliptic problem of (3.10), and (3.71), there exists a positive constant $s_3$ such that we have

$$(3.72) \qquad \begin{aligned} \|\phi^{M+1}\|_2 &\lesssim \|w^{M+1}\| + \|H^M U^{M+1}\| \\ &\lesssim \|e_w^{M+1}\| + \|w(t_{M+1})\| + \|H^M\|_{L^\infty}(\|U(t_{M+1})\| + \|e_u^{M+1}\|) \\ &\leq s_3. \end{aligned}$$

Thus, from (3.72) and (3.48), we can find a positive constant $s_4$ to get

$$(3.73) \qquad \|e_\phi^{M+1}\|_2 \leq \|\phi^{M+1}\|_2 + \|\phi(t_{M+1})\|_2 \leq s_4.$$

Furthermore, from (3.71) and (3.73), we derive

$$(3.74) \qquad \begin{aligned} \|\phi^{M+1}\|_{L^\infty} &= \|e_\phi^{M+1}\|_{L^\infty} + \|\phi(t_{M+1})\|_{L^\infty} \\ &\leq C_\Omega \|e_\phi^{M+1}\|_1^{\frac{1}{2}} \|e_\phi^{M+1}\|_2^{\frac{1}{2}} + \|\phi(t_{M+1})\|_{L^\infty} \\ &\leq C_\Omega \sqrt[4]{s_2} \sqrt{\delta t} \sqrt{s_4} + \|\phi(t_{M+1})\|_{L^\infty}. \end{aligned}$$

where we have used the Sobolev inequality $\|\phi\|_{L^\infty} \leq C_\Omega \|\phi\|_1^{\frac{1}{2}} \|\phi\|_2^{\frac{1}{2}}$ where $C_\Omega$ is a constant that only depends on $\Omega$.

Thus, if $C_\Omega \sqrt[4]{s_2} \sqrt{\delta t} \sqrt{s_4} \leq 1$, i.e., $\delta t \leq \frac{1}{C_\Omega^2 \sqrt{s_2 s_4}}$, we have

$$(3.75) \qquad \|\phi^{M+1}\|_{L^\infty} \leq 1 + \|\phi(t_{M+1})\|_{L^\infty} = \kappa.$$

Then we obtain the conclusion (3.55) by induction provided that $\delta t \leq s_0 = min(s_1, \frac{1}{C_\Omega^2 \sqrt{s_2 s_4}})$. □

**Theorem** 3.3. *Suppose the conditions of Lemma 3.4 hold, then for $0 \leq m \leq \frac{T}{\delta t} - 1$, there holds*

$$(3.76) \qquad \|e_\phi^{m+1}\|^2 + \|\nabla e_\phi^{m+1}\|^2 + \|e_u^{m+1}\|^2 + \delta t \sum_{n=0}^m (\|e_w^{n+1}\|^2 + \|\nabla e_w^{n+1}\|^2) \lesssim \delta t^2.$$

*Proof.* Since $\|\phi^n\|_{L^\infty} \leq \kappa$ for any $0 \leq n \leq \frac{T}{\delta t}$ when $\delta t \leq s_0$, by following the first step in the proof of Lemma 3.4, we obtain the conclusion (3.76). □

**4. Allen-Cahn equation.**

**4.1. Unconditional energy stable linear scheme using the IEQ approach.** For the Allen-Cahn equation, by using the same quadratization formula, we obtain a transformed PDE

system as:

$$\phi_t + M(-\epsilon^2 \Delta \phi + H(\phi)U) = 0, \tag{4.1}$$

$$U_t = \frac{1}{2} H(\phi)\phi_t, \tag{4.2}$$

where the new variable $U$ is defined as (3.1). The initial conditions $\phi|_{t=0} = \phi_0$, $U|_{t=0} = \sqrt{F(\phi_0) + B}$ and boundary conditions are (2.9). By taking the $L^2$ inner product of (4.1) with $\phi_t$, and of (4.2) with $-2U$, performing integration by parts and summing up two equalities, we can obtain the energy dissipation law of the new system (4.1)-(4.2), that reads as

$$\frac{d}{dt} E(\phi, U) = -\frac{1}{M}\|\phi_t\|^2. \tag{4.3}$$

The first-order, semi-discrete in time, IEQ scheme for solving the Allen-Cahn system (4.1)-(4.2) reads as follows,

$$\frac{\phi^{n+1} - \phi^n}{\delta t} + M(-\epsilon^2 \Delta \phi^{n+1} + H^n U^{n+1}) = 0, \tag{4.4}$$

$$U^{n+1} - U^n = \frac{1}{2} H^n (\phi^{n+1} - \phi^n), \tag{4.5}$$

where $H^n = H(\phi^n)$. The boundary conditions are:

$$(i)\ \phi^{n+1}\ \text{is periodic; or } (ii)\ \partial_{\boldsymbol{n}} \phi^{n+1}|_{\partial\Omega} = 0. \tag{4.6}$$

The unconditional energy stability of the scheme (4.4)-(4.5) is shown as follows.

**Theorem** 4.1. *The scheme* (4.4)-(4.5) *is unconditionally energy stable in the sense that*

$$E(\phi^{n+1}, U^{n+1}) \leq E(\phi^n, U^n) - \frac{1}{M\delta t}\|\phi^{n+1} - \phi^n\|^2. \tag{4.7}$$

*Proof.* By taking the $L^2$ inner product of (4.4) with $\frac{1}{M}(\phi^{n+1} - \phi^n)$ and using (3.15), we get

$$\frac{1}{M\delta t}\|\phi^{n+1} - \phi^n\|^2 + \frac{\epsilon^2}{2}(\|\nabla \phi^{n+1}\|^2 - \|\nabla \phi^n\|^2 + \|\nabla \phi^{n+1} - \nabla \phi^n\|^2) \\ + (H^n U^{n+1}, \phi^{n+1} - \phi^n) = 0. \tag{4.8}$$

By taking the $L^2$ inner product of (4.5) with $2U^{n+1}$ and using (3.15), we get

$$(\|U^{n+1}\|^2 - \|U^n\|^2 + \|U^{n+1} - U^n\|^2) = (H^n(\phi^{n+1} - \phi^n), U^{n+1}).$$

By combining the above equations together, we have

$$\frac{\epsilon^2}{2}(\|\nabla \phi^{n+1}\|^2 - \|\nabla \phi^n\|^2 + \|\nabla \phi^{n+1} - \nabla \phi^n\|^2) \\ + \|U^{n+1}\|^2 - \|U^n\|^2 + \|U^{n+1} - U^n\|^2 = -\frac{1}{M\delta t}\|\phi^{n+1} - \phi^n\|^2, \tag{4.9}$$

which concludes the energy stability (4.7) by dropping some unnecessary positive terms. □

**4.2. Implementations and well-posedness.** Similar to the Cahn-Hilliard equation, we can rewrite (4.5) as follows,

$$U^{n+1} = \frac{1}{2} H^n \phi^{n+1} + U^n - \frac{1}{2} H^n \phi^n, \tag{4.10}$$

then (4.4) can be rewritten as

$$\frac{1}{M\delta t}\phi^{n+1} - \epsilon^2 \Delta \phi^{n+1} + \frac{1}{2}H^n H^n \phi^{n+1} = \frac{1}{M\delta t}\phi^n - H^n U^n + \frac{1}{2}H^n H^n \phi^n. \tag{4.11}$$

Therefore, in practice one can solve $\phi^{n+1}$ directly from (4.11) and then update $U^{n+1}$ by (4.10).

The weak form for (4.11) can be written as the following system with unknowns $\phi \in H^1(\Omega)$,

$$\frac{1}{M\delta t}(\phi, \psi) + \epsilon^2(\nabla\phi, \nabla\psi) + \frac{1}{2}(H^n\phi, H^n\psi) = (b, \psi), \quad \psi \in H^1(\Omega). \tag{4.12}$$

where $b = \frac{1}{M\delta t}\phi^n - H^n U^n + \frac{1}{2}H^n H^n \phi^n$. We denote the above linear system as

$$(L\phi, \psi) = (b, \psi), \quad \phi, \psi \in H^1(\Omega). \tag{4.13}$$

The well-posedness of the linear system (4.13) is shown as follows.

**Theorem** 4.2. *The linear system* (4.13) *admits a unique solution* $\phi \in H^1(\Omega)$. *Furthermore, the bilinear form* $(L\phi, \psi)$ *is symmetric positive definite.*

*Proof.* (i) For any $\phi, \psi \in H^1(\Omega)$, we have

$$(L\phi, \psi) \leq \widehat{c}_1 \|\phi\|_{H^1} \|\psi\|_{H^1}, \tag{4.14}$$

where $\widehat{c}_1$ is a positive constant dependent on $\delta t, M, \epsilon$ and $\|H^n\|_{L^\infty}$. Therefore, the bilinear form $(L\phi, \psi)$ is bounded.

(ii) For any $\phi \in H^1(\Omega)$, we derive

$$(L\phi, \phi) = \frac{1}{\delta t M}\|\phi\|^2 + \epsilon^2 \|\nabla\phi\|^2 + \frac{1}{2}\|H^n\phi\|^2 \geq \widehat{c}_2 \|\phi\|_{H^1}^2, \tag{4.15}$$

□

where $\widehat{c}_2$ is a constant dependent on $\delta t, M, \epsilon$. Thus the bilinear form $(L\phi, \psi)$ is coercive. Then from the Lax-Milgram Theorem, we conclude that the linear system (4.13) admits a unique solution $\phi \in H^1(\Omega)$.

Furthermore, for any $\phi, \psi \in H^1(\Omega)$ we have $(L\phi, \psi) = (\phi, L\psi)$ that means the bilinear form is symmetric. Meanwhile, for $\phi \in H^1(\Omega)$, we have $(L\phi, \phi) \geq 0$ and the " $=$ " is valid if and only if $\phi \equiv 0$, that means the bilinear form is positive definite.

**4.3. Error estimates.** For simplicity, we still assume $\epsilon = M = 1$, and then formulate the Allen-Cahn system (4.1)-(4.2) as a truncation form:

$$\frac{\phi(t_{n+1}) - \phi(t_n)}{\delta t} - \Delta \phi(t_{n+1}) + H(\phi(t_n))U(t_{n+1}) = R_\phi^{n+1}, \tag{4.16}$$

$$U(t_{n+1}) - U(t_n) = \frac{1}{2}H(\phi(t_n))(\phi(t_{n+1}) - \phi(t_n)) + \delta t R_u^{n+1}, \tag{4.17}$$

where

$$\begin{cases} R_\phi^{n+1} = \dfrac{\phi(t_{n+1}) - \phi(t_n)}{\delta t} - \phi_t(t_{n+1}) - H(\phi(t_{n+1}))U(t_{n+1}) + H(\phi(t_n))U(t_{n+1}), \\ R_u^{n+1} = \dfrac{U(t_{n+1}) - U(t_n)}{\delta t} - U_t(t_{n+1}) + \dfrac{1}{2}H(\phi(t_{n+1}))\phi_t(t_{n+1}) - \dfrac{1}{2}H(\phi(t_n))\dfrac{\phi(t_{n+1}) - \phi(t_n)}{\delta t}. \end{cases} \tag{4.18}$$

We assume the exact solution $\phi, U$ of the system (4.1)-(4.2) possesses the following regularity conditions,

(4.19) $$\begin{cases} \phi \in L^\infty(0,T;H^2(\Omega)) \cap L^\infty(0,T;W^{1,\infty}(\Omega)), \\ U \in L^\infty(0,T;L^\infty(\Omega)), \\ \phi_t \in L^\infty(0,T;L^\infty(\Omega)), U_{tt}, \phi_{tt} \in L^2(0,T;L^2(\Omega)). \end{cases}$$

One can easily establish the following estimates for the truncation errors, provided that the exact solutions of (4.1)-(4.2) satisfy the regularity conditions (4.19).

**Lemma 4.1.** *If the exact solutions of (4.1)-(4.2) satisfy the regularity conditions (4.19), then the truncation errors satisfy*

(4.20) $$\delta t \sum_{n=0}^{[\frac{T}{\delta t}]} (\|R_\phi^{n+1}\|^2 + \|R_u^{n+1}\|^2) \lesssim \delta t^2.$$

*Proof.* Since the proof is rather standard, due to the page limit, we leave it to the interested readers. □

To derive the error estimates, we denote the error functions as

(4.21) $$e_\phi^n = \phi(t_n) - \phi^n, e_H^n = H(\phi(t_n)) - H(\phi^n), e_u^n = U(t_n) - U^n.$$

By subtracting (4.4)-(4.5) from (4.16)-(4.17), we derive the error equations:

(4.22) $$\frac{e_\phi^{n+1} - e_\phi^n}{\delta t} - \Delta e_\phi^{n+1} + e_H^n U(t_{n+1}) + H^n e_u^{n+1} = R_\phi^{n+1},$$

(4.23) $$e_u^{n+1} - e_u^n = \frac{1}{2}(e_H^n(\phi(t_{n+1}) - \phi(t_n)) + H^n(e_\phi^{n+1} - e_\phi^n)) + \delta t R_u^{n+1}.$$

Let $\kappa = \max_{0 \le t \le T} \|\phi(t)\|_{L^\infty} + 1$, we first prove the $L^\infty$ stability of solution $\phi^n$.

**Lemma 4.2.** *Suppose (i) $F(x)$ is uniformly bounded from below: $F(x) > -A$ for any $x \in (-\infty, \infty)$; (ii) $F(x) \in C^3(-\infty, \infty)$; and (iii) the exact solutions of (4.1)-(4.2) satisfy the regularity conditions (4.19), then there exists a positive constant $\widehat{s}_0$ that is given in the proof, such that when $\delta t \le \widehat{s}_0$, the solution $\phi^n$ of (4.4)-(4.5) is uniformly bounded as*

(4.24) $$\|\phi^n\|_{L^\infty} \le \kappa, n = 0, 1, \cdots, N = \frac{T}{\delta t}.$$

*Proof.* For $n = 0$, $\|\phi^0\|_{L^\infty} \le \kappa$ is true naturally. Assuming that $\|\phi^n\|_{L^\infty} \le \kappa$ is valid for all $n \le M$, we show $\|\phi^{M+1}\|_{L^\infty} \le \kappa$ is also valid through the following two steps.

(i). By taking the $L^2$ inner product of (4.22) with $e_\phi^{n+1} - e_\phi^n$, we get

(4.25) $$\frac{1}{\delta t}\|e_\phi^{n+1} - e_\phi^n\|^2 + \frac{1}{2}(\|\nabla e_\phi^{n+1}\|^2 - \|\nabla e_\phi^n\|^2 + \|\nabla e_\phi^{n+1} - \nabla e_\phi^n\|^2) \\ + (e_H^n U(t_{n+1}), e_\phi^{n+1} - e_\phi^n) + (H^n e_u^{n+1}, e_\phi^{n+1} - e_\phi^n) = (R_\phi^{n+1}, e_\phi^{n+1} - e_\phi^n).$$

By taking the $L^2$ inner product of (4.23) with $2e_u^{n+1}$, we get

(4.26) $$\|e_u^{n+1}\|^2 - \|e_u^n\|^2 + \|e_u^{n+1} - e_u^n\|^2 \\ - (e_H^n(\phi(t_{n+1}) - \phi(t_n)), e_u^{n+1}) - (H^n(e_\phi^{n+1} - e_\phi^n), e_u^{n+1}) = 2\delta t(R_u^{n+1}, e_u^{n+1}).$$

By taking the $L^2$ inner product of (4.22) with $2\delta t e_\phi^{n+1}$, we get

(4.27) $$\|e_\phi^{n+1}\|^2 - \|e_\phi^n\|^2 + \|e_\phi^{n+1} - e_\phi^n\|^2 + 2\delta t\|\nabla e_\phi^{n+1}\|^2 \\ + 2\delta t(e_H^n U(t_{n+1}), e_\phi^{n+1}) + 2\delta t(H^n e_u^{n+1}, e_\phi^{n+1}) = 2\delta t(R_\phi^{n+1}, e_\phi^{n+1}).$$

By combining the above three equalities, we derive

$$
\begin{aligned}
&\|e_\phi^{n+1}\|^2 - \|e_\phi^n\|^2 + \frac{1}{2}(\|\nabla e_\phi^{n+1}\|^2 - \|\nabla e_\phi^n\|^2) + \|e_u^{n+1}\|^2 - \|e_u^n\|^2 + 2\delta t\|\nabla e_\phi^{n+1}\|^2 \\
&+ \frac{1}{\delta t}\|e_\phi^{n+1} - e_\phi^n\|^2 + \|e_\phi^{n+1} - e_\phi^n\|^2 + \|e_u^{n+1} - e_u^n\|^2 + \frac{1}{2}\|\nabla e_\phi^{n+1} - \nabla e_\phi^n\|^2 \\
&= -(e_H^n U(t_{n+1}), e_\phi^{n+1} - e_\phi^n) + (e_H^n(\phi(t_{n+1}) - \phi(t_n)), e_u^{n+1}) \\
&\quad - 2\delta t(e_H^n U(t_{n+1}), e_\phi^{n+1}) - 2\delta t(H^n e_u^{n+1}, e_\phi^{n+1}) \\
&\quad + 2\delta t(R_\phi^{n+1}, e_\phi^{n+1}) + 2\delta t(R_u^{n+1}, e_u^{n+1}) + (R_\phi^{n+1}, e_\phi^{n+1} - e_\phi^n).
\end{aligned}
\tag{4.28}
$$

By applying Lemma 3.1, Lemma 3.2, and regularity conditions (4.19), for $n \leq M$, we estimate terms on the right hand side:

$$
\begin{aligned}
|(e_H^n U(t_{n+1}), e_\phi^{n+1} - e_\phi^n)| &\leq \|e_H^n\|\|U(t_{n+1})\|_{L^\infty}\|e_\phi^{n+1} - e_\phi^n\| \\
&\lesssim \|e_H^n\|\|e_\phi^{n+1} - e_\phi^n\| \\
&\lesssim \frac{1}{2\delta t}\|e_\phi^{n+1} - e_\phi^n\|^2 + \delta t\|e_H^n\|^2 \\
&\lesssim \frac{1}{2\delta t}\|e_\phi^{n+1} - e_\phi^n\|^2 + \delta t\|e_\phi^n\|^2;
\end{aligned}
\tag{4.29}
$$

$$
\begin{aligned}
\left|(e_H^n(\phi(t_{n+1}) - \phi(t_n)), e_u^{n+1})\right| &\leq \|e_H^n\|_{L^4}\|\phi(t_{n+1}) - \phi(t_n)\|_{L^4}\|e_u^{n+1}\| \\
&\lesssim \delta t\|e_H^n\|_{L^4}\|e_u^{n+1}\| \\
&\lesssim \delta t\|e_H^n\|_{L^4}^2 + \delta t\|e_u^{n+1}\|^2 \\
&\lesssim \delta t(\|e_H^n\|^2 + \|\nabla e_H^n\|^2) + \delta t\|e_u^{n+1}\|^2 \\
&\lesssim \delta t(\|e_\phi^n\|^2 + \|\nabla e_\phi^n\|^2) + \delta t\|e_u^{n+1}\|^2;
\end{aligned}
\tag{4.30}
$$

$$
\begin{aligned}
2\delta t|(e_H^n U(t_{n+1}), e_\phi^{n+1})| &\leq 2\delta t\|e_H^n\|\|U(t_{n+1})\|_{L^\infty}\|e_\phi^{n+1}\| \\
&\lesssim \delta t\|e_H^n\|\|e_\phi^{n+1}\| \\
&\lesssim \delta t\|e_H^n\|^2 + \delta t\|e_\phi^{n+1}\|^2 \\
&\lesssim \delta t\|e_\phi^n\|^2 + \delta t\|e_\phi^{n+1}\|^2;
\end{aligned}
\tag{4.31}
$$

$$
\begin{aligned}
2\delta t|(H^n e_u^{n+1}, e_\phi^{n+1})| &\leq 2\delta t\|H^n\|_{L^\infty}\|e_u^{n+1}\|\|e_\phi^{n+1}\| \\
&\lesssim \delta t\|e_u^{n+1}\|\|e_\phi^{n+1}\| \\
&\lesssim \delta t\|e_u^{n+1}\|^2 + \delta t\|e_\phi^{n+1}\|^2;
\end{aligned}
\tag{4.32}
$$

$$
\begin{aligned}
2\delta t|(R_\phi^{n+1}, e_\phi^{n+1}) + (R_u^{n+1}, e_u^{n+1})| &\leq 2\delta t(\|R_\phi^{n+1}\|\|e_\phi^{n+1}\| + \|R_u^{n+1}\|\|e_u^{n+1}\|) \\
&\leq \delta t(\|R_\phi^{n+1}\|^2 + \|R_u^{n+1}\|^2) + \delta t\|e_\phi^{n+1}\|^2 + \delta t\|e_u^{n+1}\|^2;
\end{aligned}
\tag{4.33}
$$

and

$$
|(R_\phi^{n+1}, e_\phi^{n+1} - e_\phi^n)| \leq \|R_\phi^{n+1}\|\|e_\phi^{n+1} - e_\phi^n\| \lesssim \frac{1}{2\delta t}\|e_\phi^{n+1} - e_\phi^n\|^2 + \delta t\|R_\phi^{n+1}\|^2.
\tag{4.34}
$$

By combining the above estimates with (4.28), we derive

$$\|e_\phi^{n+1}\|^2 - \|e_\phi^n\|^2 + \frac{1}{2}(\|\nabla e_\phi^{n+1}\|^2 - \|\nabla e_\phi^n\|^2) + \|e_u^{n+1}\|^2 - \|e_u^n\|^2 + 2\delta t\|\nabla e_\phi^{n+1}\|^2$$
$$\lesssim \delta t(\|e_\phi^n\|^2 + \|\nabla e_\phi^n\|^2 + \|e_\phi^{n+1}\|^2 + \|e_u^{n+1}\|^2) + \delta t(\|R_\phi^{n+1}\|^2 + \|R_u^{n+1}\|^2).$$

Summing up the above inequality from $n = 0$ to $m$ ($m \leq M$) and dropping some unnecessary positive terms, we get

$$\|e_\phi^{m+1}\|^2 + \frac{1}{2}\|\nabla e_\phi^{m+1}\|^2 + \|e_u^{m+1}\|^2 \lesssim \delta t \sum_{n=0}^{m}(\|e_\phi^{n+1}\|^2 + \|\nabla e_\phi^{n+1}\|^2 + \|e_u^{n+1}\|^2) + \delta t^2.$$

By Gronwall's inequality, there exist two positive constants $\widehat{s}_1, \widehat{s}_2$ such that when $\delta t \leq \widehat{s}_1$,

(4.35) $$\|e_\phi^{m+1}\|^2 + \|\nabla e_\phi^{m+1}\|^2 + \|e_u^{m+1}\|^2 \leq \widehat{s}_2 \delta t^2.$$

(ii). By using the $H^2$ regularity of elliptic problem (4.4) and the estimate (4.35), there exists a positive constant $\widehat{s}_3$, such that the following inequality holds,

(4.36)
$$\|\phi^{M+1}\|_2 \lesssim \|\frac{\phi^{M+1} - \phi^M}{\delta t}\| + \|H(\phi^M)U^{M+1}\|$$
$$\lesssim \|\frac{e_\phi^{M+1} - e_\phi^M}{\delta t}\| + \|\frac{\phi(t_{M+1}) - \phi(t_M)}{\delta t}\| + \|H(\phi^M)\|_{L^\infty}\|U^{M+1}\|$$
$$\leq \widehat{s}_3,$$

that implies, there exists a constant $\widehat{s}_4$ such that,

(4.37) $$\|e_\phi^{M+1}\|_2 \leq \|\phi(t_{M+1})\|_2 + \|\phi^{M+1}\|_2 \leq \widehat{s}_4.$$

Therefore, by (4.37) and (4.35), we obtain

(4.38)
$$\|\phi^{M+1}\|_{L^\infty} \leq \|e_\phi^{M+1}\|_{L^\infty} + \|\phi(t_{M+1})\|_{L^\infty}$$
$$\leq C_\Omega \|e_\phi^{M+1}\|_1^{\frac{1}{2}}\|e_\phi^{M+1}\|_2^{\frac{1}{2}} + \|\phi(t_{M+1})\|_{L^\infty}$$
$$\leq C_\Omega \sqrt[4]{\widehat{s}_2}\sqrt{\delta t}\sqrt{\widehat{s}_4} + \|\phi(t_{M+1})\|_{L^\infty}$$
$$\leq \kappa,$$

as long as $\delta t \leq \frac{1}{C_\Omega^2 \sqrt{\widehat{s}_2 \widehat{s}_4}}$. Thus the proof is finished by setting $\widehat{s}_0 = min(\widehat{s}_1, \frac{1}{C_\Omega^2 \sqrt{\widehat{s}_2 \widehat{s}_4}})$. □

**Theorem** 4.3. *Suppose the conditions of Lemma 4.2 hold, then for $0 \leq m \leq \frac{T}{\delta t} - 1$, there holds*

(4.39) $$\|e_\phi^{m+1}\|^2 + \|\nabla e_\phi^{m+1}\|^2 + \|e_u^{m+1}\|^2 \lesssim \delta t^2.$$

*Proof.* When $\delta t \leq \widehat{s}_0$, we have $\|\phi^n\|_{L^\infty} \leq \kappa$ for any $0 \leq n \leq \frac{T}{\delta t}$. Thus, by following the first step in the proof of Lemma 4.2, we obtain the conclusion (4.39). □

**5. Concluding Remarks.** We carry out the stability and error analysis of two first order, semi-discrete time stepping schemes for solving the Cahn-Hilliard and Allen-Cahn equations. The idea to introduce some new variables through quadratization of the IEQ method is quite different from those traditional methods like implicit, explicit, nonlinear splitting, or other various tricky Taylor expansions. In the reformulated model in terms of the new variables, all nonlinear terms can be treated semi-explicitly, which then yields a well-posed linear system that still retains the energy dissipation law. For time continuous case, the energy law in terms of the new variables is identical to that in terms of old variables, while for time discrete case, it is the first order approximation. This method is widely applicable as long as the nonlinear potential is bounded from below. The induced scheme is linear, well-posed and unconditionally energy stable. For the Allen-Cahn equation, the linear system is further symmetric positive definite.

We recall that almost all numerical schemes for various gradient flow models based on the IEQ approach or its modified version SAV approach developed by Shen et. al. in [30], had been focused on their great unconditionally energy stable properties. The quadratic formula acts as an encapsulation and thus it is difficult to obtain the two quantitative relations derived in Lemma 3.1 and 3.2, that leads to the scarce of their corresponding error estimates. In this paper, some general, sufficient conditions about the boundedness/continuity of the nonlinear functional are given in order to obtain the optimal error estimates. These conditions are naturally satisfied by the commonly use polynomial type double well potential. For the logarithmic Flory-Huggins potentials of the regularized version, these conditions are appropriate as well. By utilizing the Lipschitz property of the quadratic formula together with the mathematical inductions, we rigorously derive the optimal error estimates for the first order IEQ schemes. Moreover, the analytical approach developed in this paper is general enough and thus it can work as a standard framework to derive error estimates of IEQ type schemes for various gradient flow models with diverse nonlinear potentials. Although we consider only time discrete schemes in this paper, the results here can be carried over to any consistent finite-dimensional Galerkin type approximations since the analyses are based on the variational formulation with all test functions in the same space as the space of the trial functions.